\documentclass[english,12pt,a4paper]{article}
\usepackage[utf8]{inputenc}
\usepackage[T1]{fontenc}
\usepackage[english]{babel} 
\usepackage{csquotes}
\usepackage{charter}
\usepackage{amsmath}
\usepackage{amssymb}
\usepackage{amsthm}
\usepackage[nottoc]{tocbibind}
\usepackage{graphicx}
\usepackage{stmaryrd}
\usepackage{accents}
\usepackage[top=2.5cm, bottom=2cm, left = 2.5cm, right= 2.5cm]{geometry}
\usepackage{amsfonts}
\usepackage{dsfont}
\usepackage{authblk}
\usepackage{xcolor}
\usepackage[linesnumbered,ruled,vlined]{algorithm2e}
\usepackage{ifthen}
\usepackage{hyperref}
\usepackage[all]{hypcap}
\usepackage{aliascnt}
\usepackage[headings]{fancyhdr}
\usepackage{tabularx}
\usepackage{pgf, tikz}
\usetikzlibrary{fit,patterns,decorations.pathreplacing}
\usepackage[uniquename = false, block = none,style=authoryear,natbib=true,maxbibnames=9,maxcitenames=2,uniquelist=false,backend = biber,giveninits=true]{biblatex}
\usepackage{etoc}
\usepackage[toctitles]{titlesec}
\usepackage{caption}
\usepackage{multirow} 
\usepackage{multicol} 
\usepackage[normalem]{ulem}
\usepackage{mdframed}
\usepackage{geometry}
\usepackage{colortbl}
\usepackage{array}
\usepackage[absolute]{textpos}

\definecolor{Prune}{RGB}{99,0,60} 
\definecolor{B1}{RGB}{49,62,72} 
\definecolor{C1}{RGB}{124,135,143}
\definecolor{D1}{RGB}{213,218,223}
\definecolor{A2}{RGB}{198,11,70}
\definecolor{B2}{RGB}{237,20,91}
\definecolor{C2}{RGB}{238,52,35}
\definecolor{D2}{RGB}{243,115,32}
\definecolor{A3}{RGB}{124,42,144}
\definecolor{B3}{RGB}{125,106,175}
\definecolor{C3}{RGB}{198,103,29}
\definecolor{D3}{RGB}{254,188,24}
\definecolor{A4}{RGB}{0,78,125}
\definecolor{B4}{RGB}{14,135,201}
\definecolor{C4}{RGB}{0,148,181}
\definecolor{D4}{RGB}{70,195,210}
\definecolor{A5}{RGB}{0,128,122}
\definecolor{B5}{RGB}{64,183,105}
\definecolor{C5}{RGB}{140,198,62}
\definecolor{D5}{RGB}{213,223,61}

\DeclareNameAlias{sortname}{family-given}
\DeclareFieldFormat[article,periodical,incollection,inproceedings,proceedings,InProceedings,misc]{title}{#1}
\DeclareFieldFormat{journaltitle}{\emph{#1}}
\DeclareFieldFormat{issuetitle}{\emph{#1}}
\DeclareFieldFormat{archivePrefix}{\emph{#1}}

\renewbibmacro*{in:}{}

\usepackage[linecolor=blue!60!,backgroundcolor=blue!10!,textwidth=2.5cm,textsize=scriptsize]{todonotes}

\hypersetup{%
  pdftitle={On strongly coherent confidence envelopes for the false discoveries},            
  pdfauthor={Romain Périer, Gilles Blanchard, Guillermo Durand},          
  unicode=true,           
  colorlinks=true,       
  linkcolor = purple,
  citecolor=blue,      
}

\SetKwInput{KwInput}{Input}
\SetKwInput{KwOutput}{Output}
\SetKwInput{KwReturn}{Return}
\newcommand{\func}[5][]{
\ifthenelse{\equal {#1}{}}
{
\left\{\begin{array}{rcl}
 #2 & \longrightarrow & #3 \\
 #4 & \longmapsto & #5 \end{array}\right.
}
{#1:\left\{\begin{array}{rcl}
 #2 & \longrightarrow & #3 \\
 #4 & \longmapsto & #5 \end{array}\right.}
 }

\DeclareCiteCommand{\cite}
  {\usebibmacro{prenote}}
  {\usebibmacro{citeindex}%
   \printtext[bibhyperref]{\usebibmacro{cite}}}
  {\multicitedelim}
  {\usebibmacro{postnote}}

\DeclareCiteCommand*{\cite}
  {\usebibmacro{prenote}}
  {\usebibmacro{citeindex}%
   \printtext[bibhyperref]{\usebibmacro{citeyear}}}
  {\multicitedelim}
  {\usebibmacro{postnote}}

\DeclareCiteCommand{\parencite}[\mkbibparens]
  {\usebibmacro{prenote}}
  {\usebibmacro{citeindex}%
    \printtext[bibhyperref]{\usebibmacro{cite}}}
  {\multicitedelim}
  {\usebibmacro{postnote}}

\DeclareCiteCommand*{\parencite}[\mkbibparens]
  {\usebibmacro{prenote}}
  {\usebibmacro{citeindex}%
    \printtext[bibhyperref]{\usebibmacro{citeyear}}}
  {\multicitedelim}
  {\usebibmacro{postnote}}

\DeclareCiteCommand{\footcite}[\mkbibfootnote]
  {\usebibmacro{prenote}}
  {\usebibmacro{citeindex}%
  \printtext[bibhyperref]{ \usebibmacro{cite}}}
  {\multicitedelim}
  {\usebibmacro{postnote}}

\DeclareCiteCommand{\footcitetext}[\mkbibfootnotetext]
  {\usebibmacro{prenote}}
  {\usebibmacro{citeindex}%
   \printtext[bibhyperref]{\usebibmacro{cite}}}
  {\multicitedelim}
  {\usebibmacro{postnote}}

\DeclareCiteCommand{\textcite}
 {\boolfalse{cbx:parens}}
  {\usebibmacro{citeindex}%
    \printtext[bibhyperref]{\printnames{labelname}%
      \printtext{ (\usebibmacro{citeyear}\printtext{)}}}}
 {\ifbool{cbx:parens}
  {\bibcloseparen\global\boolfalse{cbx:parens}}
  {}%
 \multicitedelim}
{\usebibmacro{textcite:postnote}}

\let\originalleft\left
\let\originalright\right
\renewcommand{\left}{\mathopen{}\mathclose\bgroup\originalleft}
\renewcommand{\right}{\aftergroup\egroup\originalright}

\RequirePackage{xstring}

\newcommand{\paren}[2][a]{%
\IfEqCase{#1}{%
{a}{\left(#2\right)}%
{0}{(#2)}%
{1}{\big(#2\big)}%
{2}{\Big(#2\Big)}%
{3}{\bigg(#2\bigg)}%
{4}{\Bigg(#2\Bigg)}%
}[\PackageError{paren}{Undefined option to paren: #1}{}]%
}
\newcommand{\norm}[2][a]{%
\IfEqCase{#1}{%
{a}{\left\lVert#2\right\rVert}%
{0}{\lVert#2\rVert}%
{1}{\big\lVert#2\big\rVert}%
{2}{\Big\lVert#2\Big\rVert}%
{3}{\bigg\lVert#2\bigg\rVert}%
{4}{\Bigg\lVert#2\Bigg\rVert}%
}[\PackageError{norm}{Undefined option to norm: #1}{}]%
}
\newcommand{\brac}[2][a]{%
\IfEqCase{#1}{%
{a}{\left[#2\right]}%
{0}{[#2]}%
{1}{\big[#2\big]}%
{2}{\Big[#2\Big]}%
{3}{\bigg[#2\bigg]}%
{4}{\Bigg[#2\Bigg]}%
}[\PackageError{brac}{Undefined option to brac: #1}{}]%
}
\newcommand{\inner}[2][a]{%
\IfEqCase{#1}{%
{a}{\left\langle#2\right\rangle}%
{0}{\langle#2\rangle}%
{1}{\big\langle#2\big\rangle}%
{2}{\Big\langle#2\Big\rangle}%
{3}{\bigg\langle#2\bigg\rangle}%
{4}{\Bigg\langle#2\Bigg\rangle}%
}[\PackageError{inner}{Undefined option to inner: #1}{}]%
}

\newcommand{\floor}[2][a]{%
\IfEqCase{#1}{%
{a}{\left\lfloor#2\right\rfloor}%
{0}{\lfloor#2\rfloor}%
{1}{\big\lfloor#2\big\rfloor}%
{2}{\Big\lfloor#2\Big\rfloor}%
{3}{\bigg\lfloor#2\bigg\rfloor}%
{4}{\Bigg\lfloor#2\Bigg\rfloor}%
}[\PackageError{floor}{Undefined option to floor: #1}{}]%
}

\newcommand{\ceil}[2][a]{%
\IfEqCase{#1}{%
{a}{\left\lceil#2\right\rceil}%
{0}{\lceil#2\rceil}%
{1}{\big\lceil#2\big\rceil}%
{2}{\Big\lceil#2\Big\rceil}%
{3}{\bigg\lceil#2\bigg\rceil}%
{4}{\Bigg\lceil#2\Bigg\rceil}%
}[\PackageError{ceil}{Undefined option to ceil: #1}{}]%
}

\newcommand{\abs}[2][a]{
\IfEqCase{#1}{%
{a}{\left\vert#2\right\rvert}%
{0}{\vert#2\rvert}%
{1}{\big\vert#2\big\rvert}%
{2}{\Big\vert#2\Big\rvert}%
{3}{\bigg\vert#2\bigg\rvert}%
{4}{\Bigg\vert#2\Bigg\rvert}%
}[\PackageError{abs}{Undefined option to abs: #1}{}]%
}

\newcommand{\set}[2][a]{
\IfEqCase{#1}{%
{a}{\left\{#2\right\}}%
{0}{\{#2\}}%
{1}{\big\{#2\big\}}%
{2}{\Big\{#2\Big\}}%
{3}{\bigg\{#2\bigg\}}%
{4}{\Bigg\{#2\Bigg\}}%
}[\PackageError{set}{Undefined option to set: #1}{}]%
}

\newcommand{\Pro}[2][a]{\mathbb{P}\paren[#1]{#2}} 
\newcommand{\Proo}[3][a]{\mathbb{P}_{#3}\paren[#1]{#2}}

\newcommand{\Espp}[3][a]{\mathbb{E}_{#3}\brac[#1] {#2}}
\newcommand{\ind}[1]{\mathds{1}_{#1}}
\newcommand{\indev}[2][a]{\mathds{1}{\set[#1]{#2}}}
\newcommand{\EM}{\ensuremath}
\newcommand{\cA}{\EM{\mathcal{A}}}
\newcommand{\cB}{\EM{\mathcal{B}}}

\newcommand{\cH}{\EM{\mathcal{H}}}

\newcommand{\cK}{\EM{\mathcal{K}}}

\newcommand{\cP}{\EM{\mathcal{P}}}

\newcommand{\cX}{\EM{\mathcal{X}}}

\newcommand{\Nm}{\left\llbracket1,m\right\rrbracket}

\newcommand{\NN}{\left\llbracket1,N\right\rrbracket}
\newcommand{\Rfam}{\mathfrak{R}}

\newcommand{\Vstar}{V^*_{\Rfam}}

\newcommand{\RR}{\mathbb R}
\newcommand{\statfam}{\mathfrak{F}}
\newcommand{\Vhat}{\widehat{V}}

\DeclareMathOperator{\FWER}{FWER}
\DeclareMathOperator{\JER}{JER}
\DeclareMathOperator{\GHS}{GHS}
\DeclareMathOperator{\XSG}{XSG}

\DeclareMathOperator{\CT}{CT}
\DeclareMathOperator{\DBNR}{DBNR}

\def\NewTheorem#1{%
  \newaliascnt{#1}{Theoremm}
  \newtheorem{#1}[#1]{#1}
  \aliascntresetthe{#1}
  \expandafter\def\csname #1autorefname\endcsname{#1}
}
\NewTheorem{Theorem}
\NewTheorem{Definition}
\NewTheorem{Proposition-Definition}
\NewTheorem{Lemma}
\NewTheorem{Example}
\NewTheorem{Corollary}
\NewTheorem{Proposition}
\NewTheorem{Remark}

\expandafter\def\csname Modelautorefname\endcsname{Model}
\newtheorem{Assumption}{Assumption}
\expandafter\def\csname Assumptionautorefname\endcsname{Assumption}

\expandafter\def\csname Procedureautorefname\endcsname{Procedure}

\expandafter\def\csname Notationautorefname\endcsname{Notation}
\expandafter\def\csname Subsectionautorefname\endcsname{Subsection}
\addto\extrasenglish{

}
\newcommand{\secref}[1]
{\hyperref[#1]{Section~\ref{#1}}}
\newcommand{\appref}[1]
{\hyperref[#1]{Appendix~\ref{#1}}}
\newcommand{\equaref}[1]
{\hyperref[#1]{Equation~\eqref{#1}}}
\newcommand{\algoref}[1]
{\hyperref[#1]{Algorithm~\ref{#1}}}

\newcommand{\defn}[2][]{
\ifthenelse{\equal {#1}{}}
    {\begin{Definition}
        #2
    \end{Definition}}
        {\begin{Definition}[#1]
        #2
    \end{Definition}}
}

\newcommand{\defnl}[3][]{
\ifthenelse{\equal {#1}{}}
    {\begin{Definition}\label{#3}
        #2
    \end{Definition}}
        {\begin{Definition}[#1]\label{#3}
        #2
    \end{Definition}}
}

\newcommand{\thm}[2][]{
\ifthenelse{\equal {#1}{}}
    {\begin{Theorem}
        #2
    \end{Theorem}}
        {\begin{Theorem}[#1]
        #2
    \end{Theorem}}
}
\newcommand{\thml}[3][]{
\ifthenelse{\equal {#1}{}}
    {\begin{Theorem}\label{#3}
        #2
    \end{Theorem}}
        {\begin{Theorem}[#1]\label{#3}
        #2
    \end{Theorem}}
}

\newcommand{\prop}[2][]{
\ifthenelse{\equal {#1}{}}
    {\begin{Proposition}
        #2
    \end{Proposition}}
        {\begin{Proposition}[#1]
        #2
    \end{Proposition}}
}

\newcommand{\propl}[3][]{
\ifthenelse{\equal {#1}{}}
    {\begin{Proposition}\label{#3}
        #2
    \end{Proposition}}
        {\begin{Proposition}[#1]\label{#3}
        #2
    \end{Proposition}}
}

\newcommand{\rem}[2][]{
\ifthenelse{\equal {#1}{}}
    {\begin{Remark}
        #2
    \end{Remark}}
        {\begin{Remark}[#1]
        #2
    \end{Remark}}
}

\newcommand{\reml}[3][]{
\ifthenelse{\equal {#1}{}}
    {\begin{Remark}\label{#3}
        #2
    \end{Remark}}
        {\begin{Remark}[#1]\label{#3}
        #2
    \end{Remark}}
}

\newcommand{\lem}[2][]{
\ifthenelse{\equal {#1}{}}
	{
	\begin{Lemma}
	#2
	\end{Lemma}
	}
	{
	\begin{Lemma}[#1]
	#2
	\end{Lemma}
	}
}

\newcommand{\leml}[3][]{
\ifthenelse{\equal {#1}{}}
	{
	\begin{Lemma}\label{#3}
	#2
	\end{Lemma}
	}
	{
	\begin{Lemma}[#1]\label{#3}
	#2
	\end{Lemma}
	}
}

\newcommand{\cor}[2][]{
\ifthenelse{\equal {#1}{}}
	{
	\begin{Corollary}
	#2
	\end{Corollary}
	}
	{
	\begin{Corollary}[#1]
	#2
	\end{Corollary}
	}
}

\newcommand{\corl}[3][]{
\ifthenelse{\equal {#1}{}}
	{
	\begin{Corollary}\label{#3}
	#2
	\end{Corollary}
	}
	{
	\begin{Corollary}[#1]\label{#3}
	#2
	\end{Corollary}
	}
}

\newcommand{\ex}[2][]{
\ifthenelse{\equal {#1}{}}
    {\begin{Example}
        #2
    \end{Example}}
        {\begin{Example}[#1]
        #2
    \end{Example}}
}
\newcommand{\exl}[3][]{
\ifthenelse{\equal {#1}{}}
    {\begin{Example}\label{#3}
        #2
    \end{Example}}
        {\begin{Example}[#1]\label{#3}
        #2
    \end{Example}}
}

\newcommand{\ass}[2][]{
\ifthenelse{\equal {#1}{}}
    {\begin{Assumption}
        #2
    \end{Assumption}}
        {\begin{Assumption}[#1]
        #2
    \end{Assumption}}
}
\newcommand{\assl}[3][]{
\ifthenelse{\equal {#1}{}}
    {\begin{Assumption}\label{#3}
        #2
    \end{Assumption}}
        {\begin{Assumption}[#1]\label{#3}
        #2
    \end{Assumption}}
}

\newcommand{\defnfr}[2][]{
\ifthenelse{\equal {#1}{}}
    {\begin{Définition}
        #2
    \end{Définition}}
        {\begin{Définition}[#1]
        #2
    \end{Définition}}
}

\newcommand{\defnlfr}[3][]{
\ifthenelse{\equal {#1}{}}
    {\begin{Définition}\label{#3}
        #2
    \end{Définition}}
        {\begin{Définition}[#1]\label{#3}
        #2
    \end{Définition}}
}

\newcommand{\thmfr}[2][]{
\ifthenelse{\equal {#1}{}}
    {\begin{Théorème}
        #2
    \end{Théorème}}
        {\begin{Théorème}[#1]
        #2
    \end{Théorème}}
}
\newcommand{\thmlfr}[3][]{
\ifthenelse{\equal {#1}{}}
    {\begin{Théorème}\label{#3}
        #2
    \end{Théorème}}
        {\begin{Théorème}[#1]\label{#3}
        #2
    \end{Théorème}}
}

\newcommand{\propfr}[2][]{
\ifthenelse{\equal {#1}{}}
    {\begin{Proposition}
        #2
    \end{Proposition}}
        {\begin{Proposition}[#1]
        #2
    \end{Proposition}}
}

\newcommand{\proplfr}[3][]{
\ifthenelse{\equal {#1}{}}
    {\begin{Proposition}\label{#3}
        #2
    \end{Proposition}}
        {\begin{Proposition}[#1]\label{#3}
        #2
    \end{Proposition}}
}

\newcommand{\remfr}[2][]{
\ifthenelse{\equal {#1}{}}
    {\begin{Remarque}
        #2
    \end{Remarque}}
        {\begin{Remarque}[#1]
        #2
    \end{Remarque}}
}

\newcommand{\remlfr}[3][]{
\ifthenelse{\equal {#1}{}}
    {\begin{Remarque}\label{#3}
        #2
    \end{Remarque}}
        {\begin{Remarque}[#1]\label{#3}
        #2
    \end{Remarque}}
}

\newcommand{\lemfr}[2][]{
\ifthenelse{\equal {#1}{}}
	{
	\begin{Lemme}
	#2
	\end{Lemme}
	}
	{
	\begin{Lemme}[#1]
	#2
	\end{Lemme}
	}
}

\newcommand{\lemlfr}[3][]{
\ifthenelse{\equal {#1}{}}
	{
	\begin{Lemme}\label{#3}
	#2
	\end{Lemme}
	}
	{
	\begin{Lemme}[#1]\label{#3}
	#2
	\end{Lemme}
	}
}

\newcommand{\corfr}[2][]{
\ifthenelse{\equal {#1}{}}
	{
	\begin{Corollaire}
	#2
	\end{Corollaire}
	}
	{
	\begin{Corollaire}[#1]
	#2
	\end{Corollaire}
	}
}

\newcommand{\corlfr}[3][]{
\ifthenelse{\equal {#1}{}}
	{
	\begin{Corollaire}\label{#3}
	#2
	\end{Corollaire}
	}
	{
	\begin{Corollaire}[#1]\label{#3}
	#2
	\end{Corollaire}
	}
}

\newcommand{\exfr}[2][]{
\ifthenelse{\equal {#1}{}}
    {\begin{Exemple}
        #2
    \end{Exemple}}
        {\begin{Exemple}[#1]
        #2
    \end{Exemple}}
}
\newcommand{\exlfr}[3][]{
\ifthenelse{\equal {#1}{}}
    {\begin{Exemple}\label{#3}
        #2
    \end{Exemple}}
        {\begin{Exemple}[#1]\label{#3}
        #2
    \end{Exemple}}
}

\newcommand{\assfr}[2][]{
\ifthenelse{\equal {#1}{}}
    {\begin{Hypothèse}
        #2
    \end{Hypothèse}}
        {\begin{Hypothèse}[#1]
        #2
    \end{Hypothèse}}
}
\newcommand{\asslfr}[3][]{
\ifthenelse{\equal {#1}{}}
    {\begin{Hypothèse}\label{#3}
        #2
    \end{Hypothèse}}
        {\begin{Hypothèse}[#1]\label{#3}
        #2
    \end{Hypothèse}}
}

\AtEveryBibitem{
  \clearfield{doi}
  \clearfield{url}
  \clearfield{issn}
}

\title{On strongly coherent confidence envelopes for the false discoveries}

\author[1,2]{Romain Périer \thanks{Corresponding author:  \url{romain.perier@universite-paris-saclay.fr}}}
\author[1,3]{Gilles Blanchard}
\author[1,2]{Guillermo Durand}

\affil[1]{Laboratoire de Mathématiques d'Orsay (LMO), Université Paris-Saclay, CNRS}
\affil[2]{INRIA Celeste}
\affil[3]{INRIA DataShape}

\date{\today}

\begin{document}

\maketitle

Coherence is a long term notion in the literature of multiple testing. From 
$\FWER$ control multiple testing theory to confidence envelopes on the number of false positives, 
coherence is a desirable property for multiple testing procedures. However,
we show in this article that coherence is not a sufficient property to ensure that a multiple testing procedure is nonimprovable via 
an interpolation mechanism. We introduce a new notion, the notion 
of strong coherence, which is a stronger property than coherence and 
ensures that a multiple testing procedure is nonimprovable 
via the interpolation mechanisms considered in the literature.
Finally, we provide a uniform improvement of the confidence envelope from \citet{durand2020post} that is computable in polynomial time and
demonstrate the improvement on simulated data.

\localtableofcontents
\vspace{1cm}
\section{Introduction}
    \subsection{Motivations}
        Many modern statistical data analyses involve high-dimensional data, with a large number of variables 
        and many questions of interest. Frequently, scientists
        are interested in identifying a subset of variables that are responsible
        of a particular phenomenon, hence, multiple hypothesis testing theory 
        provides a natural framework to address this problem. While classical multiple hypothesis testing 
        prescribes algorithms to control specific error metrics, scientists can be tempted to snoop into data
        to find a subset of variables that are more likely to be significant, and then report the results
        of the analysis on this subset. However, this practice can lead to a significant inflation of the error metrics
        and hence to misleading conclusions. In this work, we propose a new methodology 
        relying on the one from \citet{durand2020post} to provide a valid uniform improvement of the 
        confidence bounds from the latter article. We also provide a new theoretical notion, the notion of 
        \textit{strong coherence} for confidence envelopes. A confidence envelope is strongly coherent if it is non-improvable 
        by the interpolation mechanism used in \citet{blanchard2020post} in a 
        sense defined in \eqref{def_opjer}. This notion is stronger than the notion
        of coherence from \citet{goeman2021only}, hence the name of strong coherence. 

    \subsection{Model and notations}
        For a positive integer $m$, we denote by $\Nm = \set{1,\ldots,m}$ and, for a set $A$, $\cP(A)$ the power set of $A$.
        For two real numbers $a$ and $b$, we denote by $a \wedge b$ the minimum of $a$ and $b$.\\
        We consider a statistical model $\paren{\statfam, \cX}$ where $\statfam$ 
        is a family of probability measures on a measurable space $\paren{\cX, \Xi}$. 
        We denote by $X$ an observed random variable taking values in $\cX$ defined on an abstract probability space $(\Omega, \Sigma, \mathbb{P})$,
        such that the pushforward measure of $X$ is in $\statfam$. We will denote by $X\sim \mu$ or say that $X$ follows the distribution $\mu$ 
        when the pushforward measure of $X$ is $\mu\in \statfam$. The notations $\Espp{\cdot}{\mu}$ 
        and $\Proo{\cdot}{\mu}$ denote the expectation and the probability when $X \sim \mu$.
        We consider $m \geq 1$ null hypotheses $H_1, \ldots, H_m \subseteq \statfam$. 
        For all $\mu\in \statfam$, we denote by $\cH_0(\mu) = \set{i\in \Nm : \mu\in H_i}$ 
        the set of true null hypotheses, and by $\cH_1(\mu) = \cH_0(\mu)^c$ the set 
        of false null hypotheses and by $m_0(\mu) = |\cH_0(\mu)|$ the number of true null hypotheses. 
        The dependence in $\mu$ is omitted when it is clear from the context. Further, for all $S \subseteq \Nm$, 
        we denote by $V(S) = |S\cap \cH_0|$  the number of true null hypotheses in $S$.
        With the formalism introduced up to now, we define a simultaneous confidence 
        bound for the number of false discoveries by the following functional:
        \[
            \Vhat :\left\{\begin{array}{rcl}
            \mathcal X\times (0,1) &  \longrightarrow &  \paren[2]{\cP\paren[1]{\Nm} \to \llbracket0,m\rrbracket} \\
            (x,\alpha)&\longmapsto   &  \Vhat_{(x,\alpha)}(\cdot)  \\
                \end{array},\right.
        \]
        such that the following inequality is true:
        \begin{equation}
            \forall \mu\in\statfam, \forall \alpha \in (0,1),\; \Proo{\forall S \subseteq \Nm, V(S)\leq \Vhat_{(X,\alpha)}(S)}{\mu}\geq 1-\alpha.
            \label{eq_confidence_formal2}
        \end{equation}
        The dependence in $X$ and $\alpha$ is omitted when it is clear from the context, and $\Vhat_{(X,\alpha)}$ will be written as $\Vhat$. 
        For simplicity, we assume that, for all $S \subseteq \Nm$, $\Vhat(S) \leq |S|$, as if $\Vhat$ is a valid confidence bound,
        then $\tilde{V}(\cdot) = \Vhat(\cdot) \wedge |\cdot|$ is also a valid confidence bound.

\section{Literature review} \label{sec_lit_review}
    \subsection{Interpolation devices from the \texorpdfstring{$\JER$}{JER} literature}\label{sec_interp_jer_op}
        The notion of interpolation to improve existing confidence envelopes was first introduced in \citet{blanchard2020post} and \citet{durand2020post},
        although they were applied to particular confidence envelopes, as it will be specified in \secref{sec_interp_op}. 
        In this section, we introduce these $\JER$-based interpolation devices, and we give a simpler 
        expression of the interpolation device from \citet{durand2020post}.
        Let us recall, for a reference family $\Rfam = \paren[1]{\paren{R_k(X), \zeta_k(X)}}_{k\in \cK}$ where $\cK$ is a finite set and for all
        $k\in \cK$, $R_k\subseteq \Nm$ and $0 \leq \zeta_k \leq |R_k|$, the definition of the optimal bound $\hat{V}^{\JER}_\Rfam$ from 
        \citet{blanchard2020post} and \citet{durand2020post} --- denoted $\Vstar$ in the aforementioned articles. 
        $\hat{V}^{\JER}_\Rfam$ is defined, for all $S\subseteq \Nm$, by:
        \begin{equation}\label{def_vjer}
            \hat{V}^{\JER}_\Rfam(S) = \max_{A\in \cA(\Rfam)} |S\cap A|, \text{ where } 
        \end{equation}
        \begin{equation}\label{def_A}
            \cA(\Rfam) = \set{A\subseteq \Nm : \forall k\in \cK, |A\cap R_k| \leq \zeta_k}.
        \end{equation}
        Another interpolation device was proposed by \citet{durand2020post}. This interpolation device expression 
        is given, for all $S\subseteq \Nm$, by
        \begin{equation}
            \hat{V}^{\DBNR}_\Rfam(S) = \min_{Q\subseteq\cK}\paren{\sum_{k \in Q} \zeta_k \wedge |S\cap R_k|+\abs[3]{S\setminus\bigcup_{k\in Q}R_k}}.
        \end{equation}
        We can state an equivalent but simpler expression of $\hat{V}^{\DBNR}_\Rfam$ : 
        \leml{
            Let $\Rfam = \paren[1]{\paren{R_k, \zeta_k}}_{k\in \cK}$ be a reference family. Then, for all $S\subseteq \Nm$, we have
            \begin{equation}\label{expr_dbnr}
            \hat{V}^{\DBNR}_\Rfam(S) = \min_{Q\subseteq\cK}\paren{\sum_{k \in Q} \zeta_k +\abs[3]{S\setminus\bigcup_{k\in Q}R_k}}.
            \end{equation}
        }{dbnr_simplified}
        
        We prove below expression \eqref{expr_dbnr}.
        \begin{proof}
            
            We show that there exists some $\bar{Q}\subseteq \cK$ such that the inequality $|S\cap R_k| > \zeta_k$ holds for all $k\in \bar{Q}$ and 
            $\hat{V}^{\DBNR}_\Rfam(S) = \sum_{k \in \bar{Q}} \zeta_k \wedge |S\cap R_k|+\abs[3]{S\setminus\bigcup_{k\in \bar{Q}}R_k}$.\\
            Assume that $Q\subseteq \cK$ is such that
            $\sum_{k \in Q} \zeta_k \wedge |S\cap R_k|+\abs[3]{S\setminus\bigcup_{k\in Q}R_k} = \hat{V}^{\DBNR}_\Rfam(S)$. 
            Assume that for some $\ell\in Q$, $|R_\ell \cap S| \leq \zeta_\ell$ and consider $\tilde{Q} = Q \setminus \set{\ell}$.
            The inequality
            \[
                \abs[3]{S\setminus\bigcup_{k\in Q}R_k} = \abs[3]{S\setminus\bigcup_{k\in \tilde Q}R_k} - \abs[3]{S\cap R_\ell\setminus\bigcup_{k\in \tilde Q}R_k}\geq \abs[3]{S\setminus\bigcup_{k\in \tilde Q}R_k} - \abs{S\cap R_\ell}
            \]           
            holds. Hence, we have 
            \[
                \sum_{k \in Q} \zeta_k \wedge \abs{S\cap R_k} +\abs[3]{S\setminus\bigcup_{k\in Q}R_k} \geq \sum_{k \in \tilde Q} \zeta_k \wedge \abs{S\cap R_k}+\abs[3]{S\setminus\bigcup_{k\in \tilde Q}R_k},
            \]
            Hence, $\sum_{k \in \tilde Q} \zeta_k \wedge \abs{S\cap R_k} +\abs{S\setminus\bigcup_{k\in \tilde Q}R_k} = \hat{V}^{\DBNR}_\Rfam(S)$ and, by a reasoning by induction, 
            we conclude that $\hat{V}^{\DBNR}_\Rfam(S)$ is attained for a $\bar{Q}\subseteq \cK$ such that for all $k\in \bar{Q}$, $\abs{S\cap R_k} > \zeta_k$.
            Hence, we have $\sum_{k \in \bar{Q}} \zeta_k \wedge \abs{S\cap R_k} = \sum_{k \in \bar{Q}} \zeta_k$ and 
            $\sum_{k \in \bar{Q}} \zeta_k +\abs{S\setminus\bigcup_{k\in \bar{Q}}R_k} = \hat{V}^{\DBNR}_\Rfam(S)$, which concludes the proof.

        \end{proof}
        This expression is useful in order to consider $Q\subseteq \cK$ on which $\hat{V}^{\DBNR}_\Rfam(S)$ is attained.\\
        Note that \citet{durand2020post} have shown that $\hat{V}^{\DBNR}_\Rfam \geq \hat{V}^{\JER}_\Rfam$ for any reference family $\Rfam$,
        with an equality 
        if $\Rfam$ has a forest structure, that is if for all $k, k' \in \cK$, $R_k\cap R_{k'} \in \set{\varnothing, R_k, R_{k'}}$.
        Hence, the bound $\hat{V}^{\JER}_\Rfam$ performs better than $\hat{V}^{\DBNR}_\Rfam$ on the same information basis. \\
        The idea behind these two bounds is that if we have a reference family $\Rfam$ such that its Joint Error Rate ($\JER$) is controlled, 
        that is 
        \begin{equation}\label{def_JER}
            \JER(\Rfam) =\sup_{\mu \in \statfam}\Proo{\exists k \in \cK, \abs{R_k\cap \cH_0} > \zeta_k }{X \sim \mu} \leq \alpha
        \end{equation}
        for some prespecified $\alpha \in (0,1)$, then $\hat{V}^{\JER}_\Rfam$ and $\hat{V}^{\DBNR}_\Rfam$ are valid confidence envelopes.

    \subsection{Interpolation device from \texorpdfstring{\citet{goeman2021only}}{Goeman et al. (2021)}}
        This section is dedicated to introducing the interpolation device from \citet{goeman2021only},
        the notion of \textit{coherence} defined in the aforementioned article and the notion of interpolation operator.
        We also show that the confidence envelopes $\hat{V}^{\JER}_\Rfam$ and $\hat{V}^{\DBNR}_\Rfam$ are coherent confidence envelopes. 
        \defnl{
            An interpolation operator is a deterministic 
            mapping $\mathrm{op} : \hat{V} \mapsto \mathrm{op}\paren[1]{\hat{V}}$ 
            such that if $\hat{V}$ is a confidence envelope, 
            then $\mathrm{op}\paren[1]{\hat{V}}$ is also a confidence envelope.
        }{def_interp_op}
        \citet{goeman2021only} proposed another notion of interpolation, dedicated to improve a given confidence envelope $\hat{V}$, 
        that is an operator $\mathrm{op}_{\GHS} : \hat{V} \mapsto \mathrm{op}_{\GHS}\paren[1]{\hat{V}}$ where, for 
        all $S\subseteq \Nm$, 
        \begin{equation}\label{def_opGHS}
            \mathrm{op}_{\GHS}\paren[1]{\hat{V}}(S) = \min_{U\subseteq \Nm}\set{\hat{V}(U) + \hat{V}(S\setminus U)}.
        \end{equation}
        From this notion of interpolation, \citet{goeman2021only} defines a notion of \textit{coherence} 
        for confidence envelopes that we rephrase here in terms of interpolation operators:
        \defnl{
            Let $\hat{V}$ be a confidence envelope. We say that $\hat{V}$ is coherent if $\mathrm{op}_{\GHS}\paren[1]{\hat{V}} = \hat{V}$,
            that is if for all $S\subseteq \Nm$, $\hat{V}(S) = \min_{U\subseteq \Nm}\set{\hat{V}(U) + \hat{V}(S\setminus U)}$.
        }{def_coherence}
        \reml{
            \citet{MR4771452} and \citet{goeman2021only} results are expressed in 
            a lower bound on the true discoveries fashion.
            However, it is easy to express them in term of confidence envelopes 
            on the false discoveries since, if
            $\hat d$ is a confidence envelope on the true discoveries, 
            then $\abs{\cdot} - \hat d$ is a confidence envelope on the false discoveries since, with high probability,
            for all $S\subseteq \Nm$, $\abs{S \cap \cH_0} = \abs{S} -\abs{ S \cap \cH_1} \leq \abs{S} - \hat{d}(S)$. This 
            explains that \eqref{def_opGHS} is expressed differently than in \citet{goeman2021only} and \citet{MR4771452}.
        }{rem_convert_results}

        We state below a core property on the confidence envelopes 
        $\hat{V}^{\JER}_\Rfam$ and $\hat{V}^{\DBNR}_\Rfam$ from \citet{durand2020post} and \citet{blanchard2020post}.
        \propl{
            The confidence envelopes $\hat{V}^{\JER}_\Rfam$ and $\hat{V}^{\DBNR}_\Rfam$ are coherent for any reference family $\Rfam$.
        }{coherent_envelopes}
        The proof is provided in \appref{proof_coherent_envelopes}. This property is important since it shows 
        that the interpolation device from \citet{goeman2021only} will not improve further
        the confidence envelopes $\hat{V}^{\JER}_\Rfam$ and $\hat{V}^{\DBNR}_\Rfam$. 

    \subsection{Discussion on \texorpdfstring{\citet{MR4771452}}{Xu et al. (2024)} procedure}
        This section discusses the framework from \citet{MR4771452} and its relation with the $\JER$-based interpolation devices. More
        precisely, we show that the confidence envelope $\hat{V}^{\DBNR}_\Rfam$ from 
        \citet{durand2020post} is always better than the one from \citet{MR4771452}.\\
        Further, we show that the interpolation device from \citet{goeman2021only} to define the 
        notion of coherence is not adequate because it is possible to build better confidence envelopes
        with the approach from the $\JER$ literature.\\
        From its definition of a family of $``$ focus sets$"$ $R_k$ and $``$ local bounds$"$ $\hat{V}_{R_k}$ for each $k\in \cK$ 
        --- that is a family of sets $\paren[1]{R_k, k\in \cK}$ and a family of 
        functionals $\paren[2]{\hat{V}_{R_k}, k\in \cK}$ such that for all $k\in \cK$, $\hat{V}_{R_k}$ 
        verifies \autoref{eq_confidence_formal2} with $R_k$ instead of $\Nm$ ---
        \citet{MR4771452} have defined a confidence envelope $\hat{V}^{(0)}$ by
        $\hat{V}^{(0)}\paren{R_k} = \hat{V}_{R_k}(R_k)$ for all $k\in \cK$ and
        $\hat{V}^{(0)}\paren{S} = |S|$ for all $S\subseteq \Nm$ such that $S \not \in
        \set{R_k, k\in \cK}$. Equivalently, they defined a reference family $\Rfam 
        = \paren[2]{\paren[1]{R_k, \hat{V}_{R_k}(R_k)}}_{k\in \cK}$, since no other information
        from the local bounds $\hat{V}_{R_k}, k\in \cK$ is used than the one given on the $R_k$. Hence, 
        the name $``$ focus sets$"$ is another name for the already existing reference family concept.
        Then, \citet{MR4771452} defined by induction a sequence of confidence envelopes, for all 
        $i\in \NN$, for all $S\subseteq \Nm$, by $\hat{V}^{(i+1)} (S) = \min_{U\subseteq \Nm}\set{
        \hat{V}^{(i)}(U) + \hat{V}^{(i)}(S \setminus U) }$ --- that is $\hat{V}^{(i+1)} = \mathrm{op}_{\GHS}(\hat{V}^{(i)})$ ---
        and the confidence envelope of interest is 
        $\hat{V}^{(r)}$ where $r$ is such that $\hat{V}^{(r)}=\hat{V}^{(r+1)}$ --- hence a fixed point for $\mathrm{op}_{\GHS}$, 
        that is a \textit{coherent} confidence envelope. Remark that the integer $r$ 
        is finite since the sequence is nonincreasing and bounded, and that $r$ is stochastic.\\
        From \autoref{coherent_envelopes}, we can show the following proposition,
        which states that the bound $\hat{V}^{\DBNR}_\Rfam$ is always as good as 
        $\hat{V}^{(r)}$.
        \propl{
            We have the following inequality : \[
            \hat{V}^{\JER}_\Rfam \leq \hat{V}^{\DBNR}_\Rfam \leq \hat{V}^{(r)}.\]
        }{dbnr_better}
        \begin{proof} 
            We show that $\hat{V}^{(0)} \geq \hat{V}^{\DBNR}_\Rfam$. Indeed, by definition of $\hat{V}^{(0)}$, we have
            $\hat{V}^{(0)}\paren{R_k} = \zeta_k \geq \hat{V}^{\DBNR}_\Rfam(R_k)$ for 
            all $k\in \cK$ taking $Q = \set{k}$ in \eqref{expr_dbnr},
            and $\hat{V}^{(0)}\paren{S} = |S|\geq \hat{V}^{\DBNR}_\Rfam(S)$ 
            if $S\subseteq \Nm$ is not a reference region.
            As $\hat{V}^{(0)} \geq \hat{V}^{\DBNR}_\Rfam$ and the $``$ $\geq$ $"$
            relation on confidence envelopes is stable by $\mathrm{op}_{\GHS}$, 
            we have $\hat{V}^{\DBNR}_\Rfam \leq \hat{V}^{(r)}$ because 
            $\hat{V}^{\DBNR}_\Rfam$ is coherent by \autoref{coherent_envelopes}.
        \end{proof}
        
        \reml{
            \citet{MR4771452} used the Holm-Based Combination, introduced in \citet{durand_phd}, 
            as a slight refinement of the
            Bonferroni-Based Combination, which is the tool used in \citet{durand2020post} to build the 
            reference family based on DKWM inequality. However, this improvement is over the reference family 
            --- equivalently over $\hat{V}^{(0)}$ --- but not over the interpolation operator.
        }{rem_holm_improvement}
        \reml{
            \autoref{dbnr_better} also shows that the interpolation device from \citet{goeman2021only} 
            is not optimal in the sense that with the same information 
            \textit{(i.e. the bound on the number of false positive in the reference regions, or $``$ focus sets $"$ in \citet{MR4771452} terminology)}, 
            it is possible to build a better confidence envelope from 
            \citet{durand2020post} and \citet{blanchard2020post} via the confidence envelope $\hat{V}^{\JER}_\Rfam$
            because $\hat{V}^{\JER}_\Rfam\leq \hat{V}^{\DBNR}_\Rfam$ and \citet{durand2020post} provide an example where
            $\hat{V}^{\JER}_\Rfam(S) < \hat{V}^{\DBNR}_\Rfam(S)$ for some $S$.
        }{rem_interpolation_bnr_better}
        We provide here the aforementioned example for completeness. 
        \exl{
            Let $\Nm = \set{1,2,3,4}$ and $\cK = \set{1,2,3}$ with $R_1 = \set{1,2,4}$, $R_2 = \set{2,3,4}$ and $R_3 = \set{1,3,4}$. Let $\zeta_1 = \zeta_2 = \zeta_3 =  1$. 
            Then, for $S = \Nm$, we have $\hat{V}^{\DBNR}_\Rfam(S) = 2$ and $\hat{V}^{\JER}_\Rfam(S) = 1$.
        }{ex_jer_better}
        Hence there is a gap in the notion of coherent confidence envelope. It is possible, with the same amount of information, 
        to construct a better confidence envelope with an other interpolation device than $\mathrm{op}_{\GHS}$.
        The next section is devoted to tackle this issue by introducing new interpolation operators relying
        on the $\JER$ literature. These new interpolation operators will be used to define a new notion of coherence, called \textit{strong coherence},
        which is stronger than the notion of coherence from \citet{goeman2021only}.
\section{Interpolation operators and strong coherence}
    \subsection{Interpolation operators from the \texorpdfstring{$\JER$}{JER} literature} \label{sec_interp_op}
        This section is dedicated to introduce the interpolation operators derived 
        from the confidence envelopes $\hat{V}^{\JER}_\Rfam$ and $\hat{V}^{\DBNR}_\Rfam$ from \citet{blanchard2020post} and \citet{durand2020post}.
        Let us first introduce, for a given confidence envelope $\hat{V}$, the
        reference family $\Rfam_{\hat{V}} = \paren[2]{\paren[1]{S, \hat{V}(S)}}_{S\subseteq \Nm}$ 
        that we call the reference family associated to $\hat{V}$.
        We define then two new interpolation operators, 
        $\mathrm{op}_{\JER}$ and $\mathrm{op}_{\DBNR}$, by the following 
        equations:
        \begin{equation}\label{def_opjer}
            \mathrm{op}_{\JER}\paren[1]{\hat{V}}(S) = \max_{A\in \cA(\Rfam_{\hat{V}})} \abs{S\cap A} = \hat{V}^{\JER}_{\Rfam_{\hat{V}}}(S),
        \end{equation}
        \begin{equation}\label{def_opdbnr}
            \mathrm{op}_{\DBNR}\paren[1]{\hat{V}}(S) = \min_{ \cB \subseteq \cP\paren{\Nm}}\paren{\sum_{T \in \cB} \hat{V}(T) +\abs[3]{S\setminus\bigcup_{T \in \cB}T}} = \hat{V}^{\DBNR}_{\Rfam_{\hat{V}}}(S).
        \end{equation}
        Note that the expression of $\mathrm{op}_{\DBNR}$ derives from \eqref{expr_dbnr}, and that in this special case we can again 
        simplify the expression of $\mathrm{op}_{\DBNR}$: 
        \leml{
            Let $\hat{V}$ be a confidence envelope. 
            Then, for all $S\subseteq \Nm$, we have
            \begin{equation}\label{def_opdbnr_simplified}
                \mathrm{op}_{\DBNR}\paren[1]{\hat{V}}(S) = \min_{\substack{\cB \subseteq \cP\paren{\Nm}\\ s.t. \ S \subseteq \bigcup\limits_{T \in \cB}T}}\sum_{T \in \cB} \hat{V}(T).
            \end{equation}
        }{lemma_opdbnr_simplified}

        We provide below the proof of this expression.
        \begin{proof}
            Let $\cB\subseteq \cP\paren{\Nm}$. Then
            consider $\tilde{\cB} = \cB \cup \set{S\setminus\bigcup_{T \in \cB}T}$ and we have
            \begin{align*}
                \sum_{T \in \tilde{\cB}} \hat{V}(T) +\abs[3]{S\setminus\bigcup_{T \in \tilde{\cB}}T} &\leq \sum_{T \in \cB} \hat{V}(T) + \hat{V}\paren[3]{S\setminus\bigcup_{T \in \cB}T}  \\
                &\leq \sum_{T \in \cB} \hat{V}(T) + \abs[3]{S\setminus\bigcup_{T \in \cB}T}
            \end{align*}
            where the inequality comes from the fact that $\hat{V}(T) \leq |T|$ for all 
            $T\subseteq \Nm$. Hence, we can restrict the minimization to the sets $\tilde{\cB}$ such that 
            $S\setminus\bigcup_{T \in \tilde{\cB}}T = \varnothing$, which gives the result.
        \end{proof}
        These two interpolation operators, in the same spirit as the interpolation operator 
        $\mathrm{op}_{\GHS}$ from \citet{goeman2021only},
        improve a given confidence envelope $\hat{V}$.
        It is a consequence of the following proposition.
        \propl{
            For all confidence envelope $\hat{V}$, we have 
            \[
                \mathrm{op}_{\JER}\paren[1]{\hat{V}} \leq \mathrm{op}_{\DBNR}\paren[1]{\hat{V}} \leq \mathrm{op}_{\GHS}\paren[1]{\hat{V}} \leq \hat{V}.
            \]
        }{inequality_interp_op}

        \begin{proof}
            The first inequality is a consequence of Lemma~1 from \citet{durand2020post} and the second one comes from 
            Expression~\eqref{def_opdbnr_simplified} with $\cB = \set{U,S\setminus U}$. The last inequality is a consequence of the definition
            of $\mathrm{op}_{\GHS}$ \eqref{def_opGHS} with $U = S$.
        \end{proof}
        Contrary to the interpolation operator $\mathrm{op}_{\GHS}$, the interpolation operators
        $\mathrm{op}_{\JER}$ and $\mathrm{op}_{\DBNR}$ are idempotent operators. It means
        that if we apply the interpolation operator $\mathrm{op}_{\JER}$ or $\mathrm{op}_{\DBNR}$ to
        a confidence envelope $\hat{V}$, then we cannot improve further
        the confidence envelope $\hat{V}$ by applying again the same interpolation operator. 
        The following proposition states 
        that property.
        \propl{
            The operators $\mathrm{op}_{\JER}$ and $\mathrm{op}_{\DBNR}$ are idempotent operators, that is for all confidence envelope $\hat{V}$, 
            $\mathrm{op}_{\JER}\paren{\mathrm{op}_{\JER}\paren[1]{\hat{V}}} = \mathrm{op}_{\JER}\paren[1]{\hat{V}}$ and 
            $\mathrm{op}_{\DBNR}\paren{\mathrm{op}_{\DBNR}\paren[1]{\hat{V}}} = \mathrm{op}_{\DBNR}\paren[1]{\hat{V}}$.
        }{prop_proj_opjer}
        The proof is deferred to \appref{proof_proj_opjer}.\\
        Let us define the operator $\mathrm{op}_{\XSG}$ by a rule similar to the one from \citet{MR4771452}, 
        that is, for all confidence envelope $\hat{V}$, we define $\mathrm{op}_{\XSG}\paren[1]{\hat{V}} = \lim_{r\to\infty} \mathrm{op}_{\GHS}^r\paren[1]{\hat{V}}$.
        $\mathrm{op}_{\XSG}$ exists by monotonicity of $\mathrm{op}_{\GHS}$ (\autoref{inequality_interp_op})
        and because the set of confidence envelopes is finite. 
        Then, for all confidence envelope $\hat{V}$, we have $\mathrm{op}_{\XSG}\paren[1]{\hat{V}} \geq \mathrm{op}_{\DBNR}\paren[1]{\hat{V}}$ by 
        \autoref{inequality_interp_op} which gives that $\mathrm{op}_{\DBNR}$ is nondecreasing and that $\mathrm{op}_{\DBNR}\leq \mathrm{op}_{\GHS}$, and the idempotency of $\mathrm{op}_{\DBNR}$ (\autoref{prop_proj_opjer}).

        Recall here the definition of the confidence envelope $\hat{V}^{(0)}_\Rfam$ from \citet{MR4771452} for 
        a reference family $\Rfam = \paren[2]{\paren[1]{R_k, \zeta_k}}_{k\in \cK}$, 
        that is $\hat{V}^{(0)}_\Rfam\paren{R_k} = \zeta_k$ for all $k\in \cK$ and $\hat{V}^{(0)}_\Rfam\paren{S} = |S|$
        for all $S\subseteq \Nm$ such that $S \notin \set{R_k, k\in \cK}$. \\
        An other fundamental property is that the interpolation operators $\mathrm{op}_{\JER}$ 
        and $\mathrm{op}_{\DBNR}$ are a generalization of the confidence 
        envelopes $\hat{V}^{\JER}_\Rfam$ and $\hat{V}^{\DBNR}_\Rfam$ from \citet{blanchard2020post} and \citet{durand2020post}.
        \propl{
            For all reference family $\Rfam = \paren[2]{\paren[1]{R_k, \zeta_k}}_{k\in \cK}$, we have 
            $\mathrm{op}_{\JER}\paren[1]{\hat{V}^{(0)}_\Rfam} = \hat{V}^{\JER}_\Rfam$ and $\mathrm{op}_{\DBNR}\paren[1]{\hat{V}^{(0)}_\Rfam} = \hat{V}^{\DBNR}_\Rfam$.
        }{prop_jer_dbnr_extension}
        The proof is provided in \appref{proof_jer_dbnr_extension}.

        From \autoref{prop_jer_dbnr_extension} and \autoref{ex_jer_better}, 
        we have an example of a confidence envelope $\hat{V}^{(0)}_\Rfam$ such that
        $\mathrm{op}_{\JER}\paren[1]{\hat{V}^{(0)}_\Rfam} \neq \mathrm{op}_{\DBNR}\paren[1]{\hat{V}^{(0)}_\Rfam} \leq \hat{V}^{(r)} $, 
        hence the interpolation operators $\mathrm{op}_{\JER}$ and $\mathrm{op}_{\DBNR}$ are different. Furthermore,
        it provides an example of a confidence envelope $\hat{V}^{(0)}_\Rfam$ such that
        $\mathrm{op}_{\DBNR}\paren[1]{\hat{V}}$ is a coherent confidence envelope, 
        but cannot be written as a closed testing procedure. This example is detailed in \appref{proof_ex_coherent_not_CT}.
        It implies that Corollary~1 from \citet{goeman2021only} is invalid, and that the notion of \textit{coherence}
        is not enough to define a notion of admissibility. A detailed inspection needs to be performed on the collection of claims from \citet{goeman2021only}. 
        Still, we believe that the core message of \citet{goeman2021only}, namely that only closed testing procedure are admissible, 
        remains valid with the notion of strong coherence, that we introduce in the next section, instead of the notion of coherence.\\
        Finally, we can show that the range of the interpolation operators $\mathrm{op}_{\JER}$ 
        and $\mathrm{op}_{\DBNR}$ lies in the set of coherent confidence envelopes: 
        \propl{
            We have $\mathrm{op}_{\GHS} \paren{\mathrm{op}_{\JER}} = \mathrm{op}_{\JER}$  and
            $\mathrm{op}_{\GHS} \paren{\mathrm{op}_{\DBNR}} = \mathrm{op}_{\DBNR}$.
        }{prop_proj_opjer_opdbnr}
        \begin{proof}
            The proof is a direct consequence of \autoref{coherent_envelopes} and the 
            fact that $\mathrm{op}_{\JER}\paren[1]{\hat{V}} = \hat{V}^{\JER}_{\Rfam_{\hat{V}}}$ and
            $\mathrm{op}_{\DBNR}\paren[1]{\hat{V}} = \hat{V}^{\DBNR}_{\Rfam_{\hat{V}}}$.
        \end{proof}
    \subsection{Strong coherence}
        We define here the notion of strong coherence for confidence envelopes.
        This notion is stronger than the one of \textit{coherence} from \citet{goeman2021only}, 
        and provides a correction to Corollary~1 from \citet{goeman2021only}.
        \defnl{
            Let $\hat{V}$ be a confidence envelope. 
            We say that $\hat{V}$ is strongly coherent if $\mathrm{op}_{\JER}\paren[1]{\hat{V}} = \hat{V}$, 
            that is if for all $S\subseteq \Nm$, $\hat{V}(S) = \max_{A\in \cA(\Rfam_{\hat{V}})} |S\cap A|$.
        }{def_strong_coherence}
        The following proposition shows that the notion of strong coherence is indeed 
        stronger than the notion of coherence from \citet{goeman2021only}.
        \propl{
            If a confidence envelope $\hat{V}$ is strongly coherent, 
            then it is coherent.
        }{prop_strong_coherence_implies_coherence}
        \begin{proof}
            Let $\hat{V}$ be a strongly coherent confidence envelope. Then, 
            $\hat{V} = \mathrm{op}_{\JER}\paren[1]{\hat{V}} = \hat{V}^{\JER}_{\Rfam_{\hat{V}}}$ hence 
            by \autoref{coherent_envelopes}, $\hat{V}$ is coherent. 
        \end{proof}
        Hence, every strongly coherent confidence envelope is coherent, but the converse is not true, as we 
        have already seen in \autoref{ex_jer_better}.
        Then, we show that a confidence envelope $\hat{V}$ is strongly coherent if and only if
        $\hat{V}^{\JER}_\Rfam$ for some reference family $\Rfam$.
        \propl{
            Let $\hat{V}$ be a confidence envelope. 
            Then, $\hat{V}$ is strongly coherent if and only if there exists
            a reference family $\Rfam$ such that $\hat{V} = \hat{V}^{\JER}_\Rfam$.
        }{prop_strong_coherence_as_jer}
        \begin{proof}
            Assume $\hat{V}$ is a strongly coherent confidence envelope. 
            Then, $\hat{V} = \mathrm{op}_{\JER}\paren[1]{\hat{V}} = \hat{V}^{\JER}_{\Rfam_{\hat{V}}}$ hence 
            by taking $\Rfam = \Rfam_{\hat{V}}$, we have $\hat{V} = \hat{V}^{\JER}_\Rfam$.
            Assume now that there exists a reference family $\Rfam$ such that $\hat{V} = \hat{V}^{\JER}_\Rfam$. 
            Then, \autoref{lemma_opjer_vjer} concludes that $\hat{V} = \mathrm{op}_{\JER}\paren[1]{\hat{V}}$.
        \end{proof}
        This previous proposition allows us to show that the space of closed testing procedures
        is the same as the space of strongly coherent confidence envelopes.
        We recall here the closed testing procedure for the confidence envelope on 
        false positives associated to a family of local
        tests $\mathbf{\varphi} = \paren{\varphi_S}_{S\subseteq \Nm}$, as defined in \citet{goeman2021only}.
        \defnl{
            Let $\mathbf{\varphi} = \paren{\varphi_S}_{S\subseteq \Nm}$ be a family of local tests. 
            The closed testing procedure associated to $\mathbf{\varphi}$ is the confidence envelope $\hat{V}^{\CT}_\mathbf{\varphi}$ defined by
            \begin{equation}
                \hat{V}^{\CT}_\mathbf{\varphi}(S) = \max_{A\subseteq S}\set{\abs{A}, \inf_{ U\supseteq A} \varphi_U = 0}.
            \end{equation}
        }{def_ct}
        \propl{
            Let $\hat{V}$ be a confidence envelope.
            Then, $\hat{V}$ is a closed testing procedure if and only if it is strongly coherent.
        }{prop_strong_coherence_as_CT}
        \begin{proof}
            Assume $\hat{V}$ is a closed testing procedure. Then, by \citet{blanchard2018supplement}, Section~S-1.3, there exists a reference family $\Rfam$ such that 
            $\hat{V} = \hat{V}^{\JER}_\Rfam$. By \autoref{prop_strong_coherence_as_jer}, we have that $\hat{V}$ is strongly coherent.
            Assume now that $\hat{V}$ is strongly coherent. Then, by \autoref{prop_strong_coherence_as_jer}, there exists a reference family $\Rfam$ such that 
            $\hat{V} = \hat{V}^{\JER}_\Rfam$. By \citet{blanchard2018supplement}, Section~S-1.3, we have that $\hat{V}$ is a closed testing procedure
            with the special choice of local tests $\varphi_S = 0$ if and only if $S\in \cA(\Rfam)$.
        \end{proof}
        Hence this notion of strong coherence allows us to correct Corollary~1 
        from \citet{goeman2021only} by introducing the notion of strong coherence 
        instead of the notion of coherence.
        From this last proposition, we have an alternative proof to show that the bound given in 
        \autoref{ex_jer_better} cannot be written as a closed testing procedure, since,
        in this case, $\hat{V}^{DBNR}_\Rfam$ is not strongly coherent.

\section[A uniform improvement of Durand (2020)]{Local confidence bounds: a uniform improvement of \texorpdfstring{\citet{durand2020post}}{Durand et al. (2020)} methods}
    \subsection{Considering local confidence bounds in \texorpdfstring{\citet{durand2020post}}{Durand et al. (2020)} methods}
        In this section, we show that the confidence envelope $\hat{V}^{\DBNR}_\Rfam$ from \citet{durand2020post} can be improved by
        considering, for all $k\in \cK$, a local confidence bound $\hat{V}_k$, on the reference region $R_k$ instead of only the number of false positive 
        in $R_k$ as it was done in the original work \citet{durand2020post} and 
        in \citet{MR4771452}. These two confidence envelopes, $\hat{V}^{\DBNR}_\Rfam$ and $\hat{V}^{(r)}$ 
        only use the upper bound on the number of false positive $V(R_k)$ in the reference regions
        given by the local confidence bounds. We can improve upon this by using the local 
        confidence bounds, with $S\subseteq R_k$, to bound $ V(S\cap R_k)$.\\
        The main idea is that for top-$k$ and $k$-FWER confidence bounds methods, following the analysis from \citet{périer2026confidenceenvelopesfalsediscoveries}
        concerning adaptive top-$k$ and $k$-FWER methods, 
        the event used to compute the number of false positive in $R_k$ --- that is $\zeta_k = \hat{V}_k(R_k)$ --- provides, for free, 
        a local top-$k$ structure within $R_k$. 
        Hence, instead of only considering the number of false positive in $R_k$, 
        we can consider the local confidence bounds $\hat{V}_k(S\cap R_k)$ for all $S\subseteq R_k$.\\
        Moreover, the computation of those local confidence envelopes usually requires the same computational cost than
        the cost to compute $\zeta_k$,
        since it is mandatory to prior compute $\zeta_k$ to get $\hat{V}_k(S\cap R_k)$ for some $S\subseteq R_k$. \\
        The following proposition shows that if, for all $k\in \cK$, we have local 
        confidence bounds on the reference regions $R_k$ at level $\alpha_k$,
        such that 
        $\sum\limits_{k\in \cK} \alpha_k$ is upper bounded by $\alpha$, 
        the informal construction above provides a confidence envelope $\tilde{V}$ 
        that is an $\alpha$-confidence envelope.
        \propl{
            Let $\paren{R_k}_{k\in\cK}$ be reference regions,
            let $\paren{\alpha_k}_{k\in \cK}$ such that $\sum\limits_{k\in \cK} \alpha_k \leq \alpha$,
            and for all $k\in \cK$, let $\hat{V}_k$ be a local $\alpha_k$-confidence envelope over $R_k$, 
            meaning that with probability $1 -\alpha_k$, for all $S\subseteq R_k$, $V(S) \leq \hat{V}_k(S)$.\\
            Then the bound $\tilde{V}$ defined by 
            \begin{equation}\label{def_tildeV}
                \tilde{V}(S) = \min_{Q \subseteq \cK } \paren{ \sum_{k\in Q}\hat{V}_k \paren{S \cap R_k} + \abs[3]{S\setminus \bigcup_{k\in Q} R_k}}
            \end{equation}
            is an $\alpha$-confidence envelope. \\
            Further, if for all $k\in \cK$, $\hat{V}_k(S \cap R_k) \leq \hat{V}_k(R_k)$, 
            then $\tilde{V}$ is a uniform improvement of the confidence envelope $\hat{V}^{\DBNR}_\Rfam$, where 
            $\Rfam = \paren[2]{\paren[1]{R_k, \hat{V}_k(R_k)}}_{k\in \cK}$.
        }{DBNRrefinement}
        \begin{proof}
            \begin{itemize}
            \item Let $A = \bigcup\limits_{k\in \cK} \set{\exists S \subseteq R_k : V(S) > \hat{V}_k \paren{S \cap R_k} }$. 
            By union bound, $\Pro{A} \leq \sum \limits_{k\in \cK} \alpha_k \leq \alpha$.
            \item Within the event $A^c$, for all $Q \subseteq \cK$, 
            \begin{align*}
                V(S) & \leq \sum\limits_{k\in Q} \abs{R_k \cap S \cap \cH_0 } + \abs[3]{\paren{S\cap \cH_0} \setminus \bigcup_{k\in Q} R_k}\\
                &\leq \sum\limits_{k\in Q} \hat{V}_k\paren{S \cap R_k} + \abs[3]{S\setminus \bigcup_{k\in Q} R_k}.
            \end{align*}
            The first part of this proposition follows by considering the minimum over this inequality.
            \end{itemize}
            The second part of the proposition comes by comparing the definition of $\tilde{V}$ and $\hat{V}^{\DBNR}_\Rfam$.
        \end{proof}
        This confidence envelope
        is naturally a uniform improvement of the confidence envelope $\hat{V}^{\DBNR}_\Rfam$ 
        from \citet{durand2020post} since for all $k\in \cK$, $\hat{V}_k(R_k) \geq \hat{V}_k(R_k\cap S)$ 
        for all $S\subseteq \Nm$ when the local bounds are coherent by Lemma~3 from \citet{goeman2021only} --- as are the top-$k$ and
        $k$-$\FWER$ methods described in \citet{périer2026confidenceenvelopesfalsediscoveries} 
        since they are strongly coherent hence coherent by \autoref{prop_strong_coherence_implies_coherence}.
        Thus considering these underlying structures will always be better than
        only bounding the number of false positive in $R_k$ by $\hat{V}_k(R_k)$.

        The following proposition shows that, when the local confidence envelopes 
        are provided by the JER constructions, the confidence envelope $\tilde{V}$ from 
        \autoref{DBNRrefinement} is lower bounded by the strongly coherent associated JER envelope. However, 
        this strongly coherent JER envelope is not computable in general, while $\tilde{V}$ is computable
        if the reference regions have a forest structure, see Algorithm~1 from \citet{durand2020post}.

        \prop{
            Let $\paren{R_k}_{k\in \cK}$ be reference regions,
            and, 
            for all $k\in \cK$, $\Rfam_k = \paren{R_{k,i},\zeta_{k,i}}_{i \in \cK_k}$ 
            such that for all $i \in \cK_k$, $R_{k,i} \subseteq R_k$. 
            Then, defining $\tilde{\cK} = \bigcup\limits_{k\in \cK} \set{\paren{k,i}, i\in \cK_k}$ 
            and $\tilde{\Rfam} = \paren{R_{k,i},\zeta_{k,i}}_{(k,i) \in  \tilde{\cK}}$, 
            \begin{equation}\label{ineq_jer_local}
                \hat{V}^{\JER}_{\tilde{\Rfam}}(S) \leq  \min_{Q \subseteq \cK} \paren{ \sum_{k\in Q}\hat{V}^{\JER}_{\Rfam_k}\paren{S \cap R_k} 
                + \abs[3]{S\setminus \bigcup_{k\in Q} R_k}}.
            \end{equation}
        }
        \begin{proof}
            Let $A \in \cA(\tilde{\Rfam}) = \set{A\subseteq \Nm : \forall (k,i) \in \tilde{\cK}, \abs{A\cap R_{k,i}}\leq \zeta_{k,i}}$ 
            and $Q \subseteq \cK$.
            Then, for all $k\in \cK$, 
            $A\cap R_k \in \set{B\subseteq R_k : \forall i \in \tilde{\cK}, \abs{A\cap R_{k,i}}\leq \zeta_{k,i}}$,
            thus 
            \begin{align*}
                \abs{A\cap S} &= \sum\limits_{k\in Q} \abs{ S \cap R_k \cap A } + \abs[3]{S\cap A\setminus \bigcup_{k\in Q} R_k} \\
                &\leq \sum\limits_{k\in Q} \hat{V}^{\JER}_{\Rfam_k}\paren{S \cap R_k} + \abs[3]{S\setminus \bigcup_{k\in Q} R_k}.
            \end{align*}
            It implies $ \hat{V}^{\JER}_{\tilde{\Rfam}} \paren{S}\leq \sum\limits_{k\in Q} \hat{V}^{\JER}_{\Rfam_k}\paren{S \cap R_k} + \abs[3]{S\setminus \bigcup_{k\in Q} R_k}$ and the inequality follows taking the minimum over $Q$.
        \end{proof}
        \reml{
            As it is the case for $\hat{V}^{\DBNR}_\Rfam$, if the reference 
            regions $\paren{R_k}_{k\in \cK}$ have a complete forest structure --- 
            that is the union of the regions equals $\Nm$ and that each node from the forest structure 
            is either a leaf or partitioned by its children (\textit{i.e.,} the sets $R_\ell \subsetneq R_k $) --- 
            then the minimum over $Q\subseteq \cK$ in \eqref{def_tildeV} and \eqref{ineq_jer_local}
            can be restricted to a minimum over the sets $Q\subseteq \cK$ 
            such that $\paren{R_k}_{k\in Q}$ is a partition of $\Nm$.
            Practically, for a specific set $S \subseteq \Nm$, we are applying the Algorithm~1 from \citet{durand2020post} 
            with $\hat{V}_k\paren[1]{S\cap R_k}$ instead of $\zeta_k$ for all 
            $k\in \cK$.
        }{rem_complete_forest}
    \subsection{Simulations}
        In this section, we provide simulations to illustrate the improvement of the confidence envelope $\tilde{V}$ 
        from \autoref{DBNRrefinement} over the confidence envelope $\hat{V}^{\DBNR}_\Rfam$ from \citet{durand2020post} 
        --- and the confidence envelope from \citet{MR4771452}, by extension. 
        The code to reproduce the simulations is available at \url{https://github.com/RomainPerier/Local_CB}.

        We consider the same simulation setting than in \citet{périer2026confidenceenvelopesfalsediscoveries}.
        This setting is an \textit{heterogeneous} setting, meaning that the $p$-values are not identically distributed,
        that we recall here for completeness.

        We simulate two groups of $N = 50$ subjects, and for each subject, we observe $m$ independent binary responses (\textit{e.g.} treatment effects, gene variants, ...). 
        We denote by $\beta_{1i}$ and $\beta_{2i}$ the success probability, respectively 
        for the first group and the second group. Then, we simultaneously test the null 
        hypotheses $H_i : $``$ \beta_{1i}=\beta_{2i}"$ by
        two-sided Fisher's exact tests --- see \citet{MR3087436} for more details on these tests --- 
        producing for each $i \in \Nm$ a $p$-value $p_i$ with its \textit{c.d.f.} denoted $F_i$. 
        We will then have three distinct sets of $p$-values. 
        The first set of 
        $p$-values is generated by underlying data such that 
        $\beta_{1i}=\beta_{2i} = 0.01$. Hence it is a set of $p$-values such 
        that the underlying null hypotheses are true.
        The second set of 
        $p$-values is generated 
        by underlying data such that $\beta_{1i}=\beta_{2i} = 0.1$ hence this is a set of $p$-values such
        that the underlying null hypotheses are true.
        Each of these two sets of $p$-values is equally represented in the simulations,
        meaning that $50\%$ of the $p$-values under the null hypothesis belong to the first set and $50\%$ 
        to the second set. 
        Finally, the third set of $p$-values is generated 
        by underlying data such that $\beta_{1i}= 0.1$ and $\beta_{2i} = q$, 
        where $q \in \set{0.2, 0.3, 0.4}$ represents the strength of the signal. 

\begin{figure}[h!]
            \begin{center}
            \begin{tikzpicture}[scale=0.85]
            \tikzstyle{quadri}=[circle,draw,text=black, thick]
            \tikzstyle{quadriH1}=[circle,draw,text=black,fill=pink, thick]
            \tikzstyle{estun}=[<-,>=latex,very thick]

            \node[quadri] (R18) at (0,6) {$R_1$};

            \node[quadri] (R14) at (-4,4) {$R_2$};
            \node[quadri] (R58) at (4,4) {$R_3$};

            \node[quadriH1] (R12) at (-6,2) {$R_4$};
            \node[quadri] (R34) at (-2,2) {$R_5$};
            \node[quadri] (R56) at (2,2) {$R_6$};
            \node[quadri] (R78) at (6,2) {$R_7$};

            \node[quadriH1] (R1) at (-7,0) {$R_8$};
            \node[quadriH1] (R2) at (-5,0) {$R_9$};
            \node[quadri] (R3) at (-3,0) {$R_{10}$};
            \node[quadri] (R4) at (-1,0) {$R_{11}$};
            \node[quadri] (R5) at (1,0) {$R_{12}$};
            \node[quadri] (R6) at (3,0) {$R_{13}$};
            \node[quadri] (R7) at (5,0) {$R_{14}$};
            \node[quadri] (R8) at (7,0) {$R_{15}$};

          

            \draw[estun] (R1)--(R12);
            \draw[estun] (R2)--(R12);
            \draw[estun] (R3)--(R34);
            \draw[estun] (R4)--(R34);
            \draw[estun] (R5)--(R56);
            \draw[estun] (R6)--(R56);
            \draw[estun] (R7)--(R78);
            \draw[estun] (R8)--(R78);

            \draw[estun] (R12)--(R14);
            \draw[estun] (R34)--(R14);
            \draw[estun] (R56)--(R58);
            \draw[estun] (R78)--(R58);

            \draw[estun] (R14)--(R18);
            \draw[estun] (R58)--(R18);

            \end{tikzpicture}
            \end{center}
            \caption[Perfect binary tree structures]{Perfect binary tree structures, here with a number of active leaves equal to 2. 
            The pink nodes are those containing signal. 
            This figure is one of the figures from \citet{durand2020post}, 
            but slightly modified.}
            \label{fig_binary_tree}
        \end{figure}
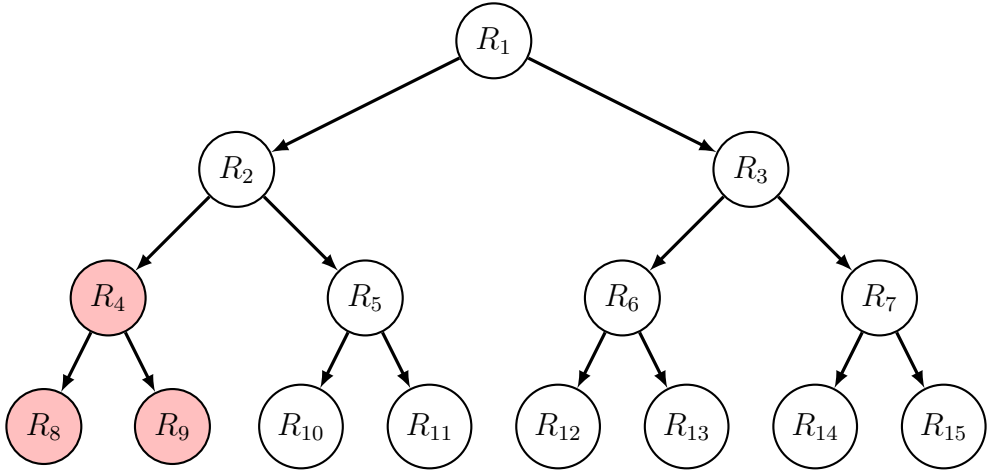

        These $p$-values are localized in a binary tree structured reference family $\Rfam$. We give 
        an example of such a structure in \autoref{fig_binary_tree}, albeit simpler than the 
        one we used. 
        For each of the 
        leaves of the tree --- that is, for each reference region $R_k$, $k\in \cK$, such that 
        $ R_k \cap R_{k'} \in \set{\varnothing, R_k}$ for all $k'\in \cK$ --- $d_k \in 
        \set{0,0.2, 0.5, 0.8}$ represents the proportion of $p$-values under the alternative hypothesis within
        the leaf. This proportion is common to all active leaves, that is,
        there exists a common value $d\in \set{0.2, 0.5, 0.8}$ such that
        for any leaf $R_k$, $k\in \cK$, verifying $d_k \neq 0$, $d_k = d$.

        Finally, for all of our 
        simulations, the proportion $d$ and the strength of the signal 
        $q$ are specified and we set the number of hypotheses to
        $ m = 3200$. In addition,
        the number of active leaves is set to $8$ and the number of hypotheses within each leaf is set to $100$, therefore the 
        forest structure possesses $32$ leaves.

        Furthermore, for each couple $(d,q)$ we draw $m$ 
        contingency tables, compute the associated $p$-values
        and apply four methods to this unique set of $p$-values. 
        For each of the $p$-values $p_i$, the 
        cumulative distribution function \textit{c.d.f} is 
        kwown and denoted $F_i$.
        We describe below the four methods that we consider in our simulations, $\hat{V}^{(1)}, \hat{V}^{(2)}, \hat{V}^{(3)}, \hat{V}^{(4)}$, but before 
        we introduce some notations.
        For a reference family $\Rfam = \paren[2]{\paren{R_k,\zeta_k}}_{k\in \cK}$, 
        we denote $K = |\cK|$ and $m_k = |R_k|$ for all $k\in \cK$.
        For a vector $x\in \RR^m$ and $A\subseteq \Nm$, 
        $x_{[i:A]}$ denotes the $i$-th order statistic of the vector $\paren[1]{x_j}_{j\in A}$ in descending order and 
        $x_{(i:A)}$ denotes the $i$-th order statistic of the vector $\paren[1]{x_j}_{j\in A}$ in ascending order.
        The methods $\hat{V}^{(1)}, \hat{V}^{(2)}$ are derived from the DKW inequality \citep{massart1990tight,dvoretzky1956asymptotic};
        $\hat{V}^{(1)}$ is $\hat{V}^{\DBNR}_\Rfam$ where $\Rfam = \paren[2]{\paren[1]{R_k, \zeta_k}}_{k\in \cK}$ with $\zeta_k $ defined by 
        \[
        \zeta_k=\min_{t\in[0,1)}\floor[a]{\paren[a]{\frac{\lambda_{\alpha/K}}{2(1-t)}+\sqrt{\frac{\lambda_{\alpha/K}^2}{4(1-t)^2} +  \frac{%
        |R_k|-i_{R_k}(t)  }{1-t}}}^2   },
        \]
        where $\lambda_{\alpha/K}=\sqrt{\frac12\log \frac K\alpha}$. This overestimator on the number of false discoveries in $R_k$ has been considered for instance 
        in \citet{durand2020post}, \citet{MBR2024} and \citet{périer2026confidenceenvelopesfalsediscoveries}. 
        $\hat{V}^{(2)}$ is defined by \eqref{def_tildeV} with, 
        for all $k\in \cK$, $\hat{V}_k$ defined
        by $\hat{V}_k =\hat{V}_{\Rfam_k}^{\JER}$ where $\Rfam_k = \paren[2]{\paren[1]{R_{k,i}, \zeta_{k,i}}}_{i\in \llbracket 1, m_k\rrbracket}$ 
        is such that for all $i\in \llbracket 1, m_k\rrbracket$, $R_{k,i} = \set{j\in R_k : p_j \leq p_{(i:R_k)} }$ 
        and $\zeta_{k,i} = \left\lfloor  \zeta_k p_{(i:R_k)}+ \sqrt{\zeta_k}\lambda_{\alpha/|\cK|}\right\rfloor$ --- hence the adaptive 
        DKW procedure described in \citet{périer2026confidenceenvelopesfalsediscoveries}, Proposition~4.14,
        but restricted to the set $R_k$.
        
        The last two methods, $\hat{V}^{(3)}$ and $\hat{V}^{(4)}$, are methods derived from the Bretagnolle inequality \citep{shorack2009empirical};
        $\hat{V}^{(3)}$ is $\hat{V}^{\DBNR}_{\bar{\Rfam}}$ where $\bar{\Rfam} = \paren[1]{\paren[0]{R_k, \bar{\zeta}_k}}_{k\in \cK}$ with $\bar{\zeta}_k $ defined by
        \[
            \bar{\zeta}_k = \max\set{n \in \llbracket 0, m_k\rrbracket : \forall t \in [0,1],  
            \sum\limits_{1\leq j\leq n}\paren[1]{\mathbb{G}(t)}_{(j: R_k)}\leq
            \lambda_{\alpha/ (eK)}\sqrt{n}},
        \]
        where $e = \exp(1)$
        and for all $t\in [0,1]$, 
        $\mathbb{G}(t) = \paren[1]{\indev[a]{p_i\leq t}- F_i(t)}_{1\leq i \leq m}$. 
        This overestimator on the number of false discoveries in $R_k$ has been considered 
        in \citet{périer2026confidenceenvelopesfalsediscoveries}.\\
        $\hat{V}^{(3)}$ is defined by \eqref{def_tildeV} with, for all $k\in \cK$, $\hat{V}_k = \bar{V}_k$ defined by
        $\bar{V}_k =\hat{V}_{\bar{\Rfam}_k}^{\JER}$ where $\bar{\Rfam}_k = \paren[2]{\paren[1]{R_{k,i}, \bar{\zeta}_{k,i}}}_{i\in \llbracket 1, m_k\rrbracket}$
        $\bar{\zeta}_{k,i}= \bar{\zeta}_k \wedge \left\lfloor \sum _{1 \le i \le \bar{\zeta}_k}
        \paren[1]{\mathbf{F}\paren[0]{p_{(k)}}}_{[i : R_k]}+ \sqrt{\bar{\zeta}_k}
        \lambda_{\alpha/ (eK)}\right\rfloor.$
        where for all $t\in [0,1]$, $\mathbf{F}(t) =  \paren[1]{F_i(t)}_{1\leq i \leq m}$.
        Hence, for all $k\in \cK$, the procedure $\bar{V}_k$ is the adaptive top-$k$ procedure from  
        \citet{périer2026confidenceenvelopesfalsediscoveries}, Theorem~6.2, but restricted to 
        $R_k$.

        \begin{figure}[ht]
        \begin{center}
        \includegraphics[scale=0.50]{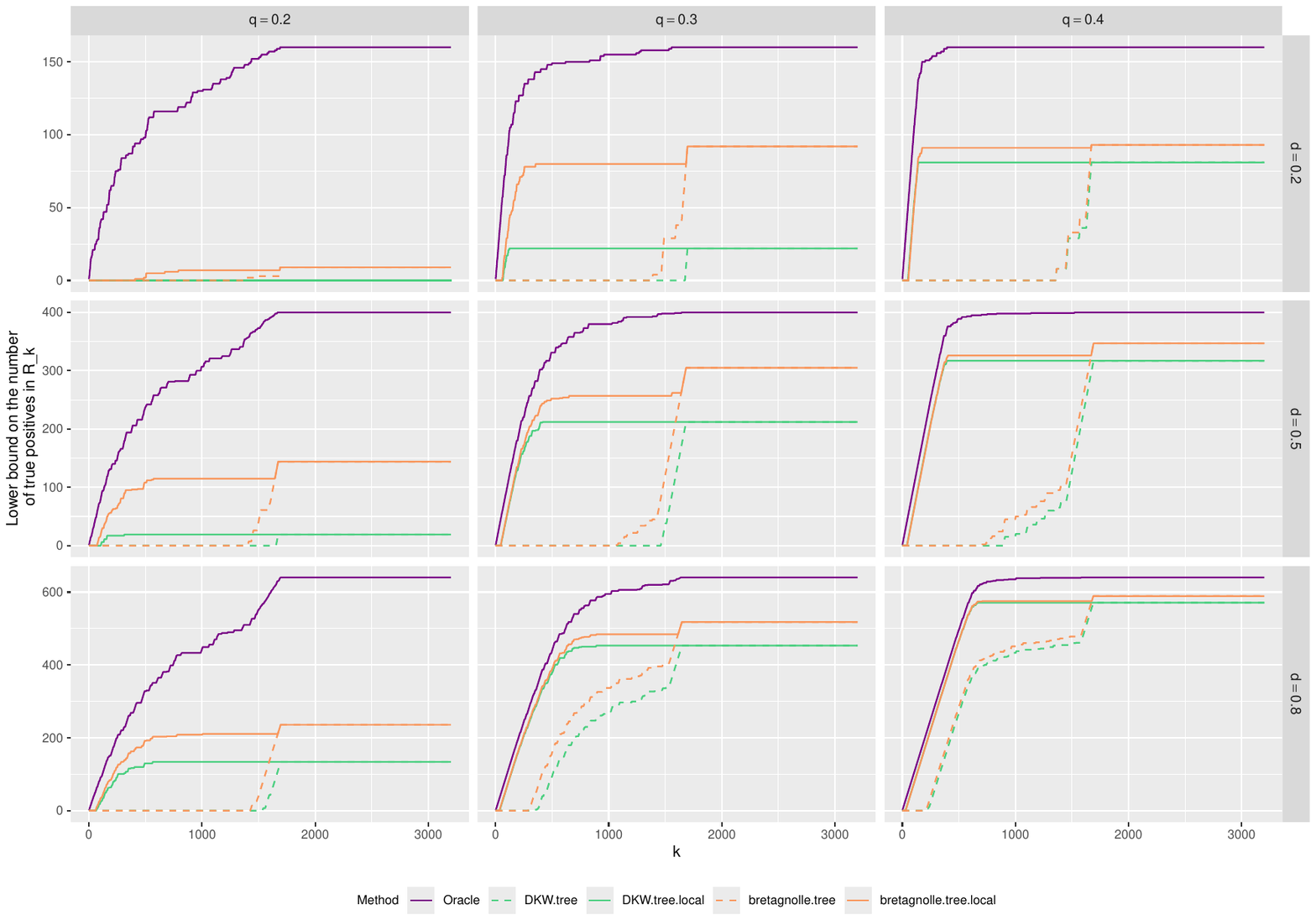}
        \end{center}
        \caption{Lower bounds on number of true positives in the top-$k$ sets provided by DKW based confidence envelopes. Higher is better.}
        \label{LBlocal}
        \end{figure}

        The results of these simulations are summerized
        by \autoref{LBlocal}.
        where we plot the lower bounds on the number of true positives for the top-$k$ sets $S_k$ with $k\in \Nm$, that is, the 
        hypotheses corresponding to the $k$ smallest $p$-values.
        Therefore, in \autoref{LBlocal} at abscissa $k$, we plot $|S_k| - \hat{V}^{(i)}(S_k)$, for all
        $i \in \set{1,2,3,4}$.
        Each of the procedures based on $\hat{V}^{\DBNR}$ are highly outperformed
        by the corresponding procedures based on $\tilde{V}$, which is consistent with \autoref{DBNRrefinement}.
        In particular, a caveat of the $\hat{V}^{\DBNR}$ procedures is that it performs 
        poorly when density of the signal is not high enough, like it is the case for our 
        simulations, because the smallest $p$-values may not be in the same reference region.

\section{Conclusion}
        In this article, we have shown that the notion of coherence from \citet{goeman2021only} is not sufficient 
        ensure that a confidence envelope can be written as a closed testing procedure. 
        We have introduced the notion of strong coherence, which is a stronger notion than coherence, and we have shown that
        a confidence envelope is strongly coherent if and only if it can be written as a JER
        confidence envelope. We have also provided a uniform improvement of the confidence envelope from \citet{durand2020post}, and thus from \citet{MR4771452},
        by considering local confidence bounds on the reference regions. Although this improvement is not equal to the optimal JER confidence envelope, 
        it is computable in practice and provides a significant improvement over the original methods, as illustrated by our simulations, using methods
        derived from DKW and Bretagnolle inequalities.
        At some point, it would be interesting to investigate other inequalities to derive confidence envelopes based on local 
        confidence bounds, and to compare their performance with the methods used without structure, as the confidence envelopes from \citet{périer2026confidenceenvelopesfalsediscoveries}. 

\section*{Acknowledgements}
    The authors acknowledge funding from the French research grants ANR-19-CHIA-0021 (BISCOTTE), ANR-21-CE23-0035 (ASCAI) and ANR-23-CE40-0018 (BACKUP).
    Romain Périer and Guillermo Durand furthermore acknowledge funding from the French research grant ANR-20-IDEES-0002 (PIA).
    A similar notion of strong coherence was independently introduced by \citet{tran2026interpolationcoherencemultipletesting}, called perfect coherence. The present 
    article is part of the PhD thesis of Romain Périer \citep{Perier2026},
    and has been developed independently from the 
    work of \citet{tran2026interpolationcoherencemultipletesting}.
    The authors declare no financial nor commercial conflicts of interests and certify 
    that none of the proofs, ideas or results has been generated 
    by generative artificial intelligence.

\printbibliography[] 

\appendix

\section{Deferred proofs}
    \subsection{Proof of \autoref{prop_proj_opjer}} \label{proof_proj_opjer}
        We first state a lemma that will be useful in the proof of \autoref{prop_proj_opjer}.
        \leml{
            Let $\Rfam = \paren{R_k,\zeta_k}_{k\in \cK}$ be a reference family and $\hat{V}^{\JER}_\Rfam$ the confidence envelope associated to $\Rfam$. 
            Then we have $\mathrm{op}_{\JER}\paren[1]{\hat{V}^{\JER}_\Rfam} = \hat{V}^{\JER}_\Rfam$.
        }{lemma_opjer_vjer}
        \begin{proof}
            We show that $\cA(\Rfam) = \cA\paren{\Rfam_{\hat{V}_\Rfam^{\JER}}}$.
            The key property that we use here comes from Proposition~2.1 from 
            \citet{blanchard2020post}: for any $A\subseteq \Nm$, $A\in \cA(\Rfam)$ 
            if and only if for all $S\subseteq \Nm$, $\abs{A\cap S} \leq \hat{V}^{\JER}_\Rfam(S)$.
            Hence, we have $\cA(\Rfam) = \cA\paren{\Rfam_{\hat{V}_\Rfam^{\JER}}}$ 
            since $\cA\paren{\Rfam_{\hat{V}_\Rfam^{\JER}}}$ is defined 
            as the set of all $A\subseteq \Nm$ such that for all $S\subseteq \Nm$, $\abs{A\cap S} \leq \hat{V}^{\JER}_\Rfam(S)$.
        \end{proof}
        
        Let $\hat{V}$ be a confidence envelope and $S\subseteq \Nm$. Denote by $\Rfam_{\hat{V}} = \paren{S, \hat{V}(S)}_{S\subseteq \Nm}$ the 
        reference family associated to $\hat{V}$.
        We first show that $\mathrm{op}_{\JER}$ is an idempotent operator. 
        By \autoref{lemma_opjer_vjer}, 
        we have, since $\hat{V}^{\JER}_{\Rfam_{\hat{V}}} = \mathrm{op}_{\JER}\paren[1]{\hat{V}}$ by definition of $\mathrm{op}_{\JER}$ in \eqref{def_opjer}, that
        $\mathrm{op}_{\JER}\paren[2]{\mathrm{op}_{\JER}\paren[1]{\hat{V}}} = 
        \mathrm{op}_{\JER}\paren[1]{\hat{V}}$.
    
        We then show that $\mathrm{op}_{\DBNR}$ is an idempotent operator. We denote 
        by $\Rfam_{\mathrm{op}_{\DBNR}\paren[1]{\hat{V}}} = \paren[2]{\paren[1]{S, \mathrm{op}_{\DBNR}\paren[1]{\hat{V}}(S)}}_{S\subseteq \Nm}$ 
        the reference family associated to $\mathrm{op}_{\DBNR}\paren[1]{\hat{V}}$. Then, for $S\subseteq \Nm$, we have
        \begin{itemize}
            \item $\hat{V}^{\DBNR}_{\Rfam_{\mathrm{op}_{\DBNR}\paren{\hat{V}}}}(S)\leq \hat{V}^{\DBNR}_{\Rfam_{\hat{V}}}(S)$ by 
            definition of $\mathrm{op}_{\DBNR}$ in \eqref{def_opdbnr} with $\cB = \set{S}$.
            \item Let $\cB \subseteq \cP(\Nm)$ such that 
            $\hat{V}^{\DBNR}_{\Rfam_{\mathrm{op}_{\DBNR}\paren{\hat{V}}}}(S) = \sum_{T\in B} \mathrm{op}_{\DBNR}\paren[1]{\hat{V}}(T)$
            and for all $T\in \cB$, 
            let $\cB_T \subseteq \cP(\Nm)$ such that 
            $\mathrm{op}_{\DBNR}\paren[1]{\hat{V}}(T) = \sum_{U\in \cB_T} \hat{V}(U)$.
            Recall that, by \eqref{def_opdbnr_simplified}, such sets exist and we have 
            that for all $T\in \cB$, $T \subseteq \bigcup\limits_{U\in \cB_T} U$.
            Then
            \begin{align*}
                \hat{V}^{\DBNR}_{\Rfam_{\mathrm{op}_{\DBNR}\paren{\hat{V}}}}(S) &= \sum_{T\in \cB} \mathrm{op}_{\DBNR}\paren[1]{\hat{V}}(T)\\
                & = \sum_{T\in \cB} \sum_{U\in \cB_T} \hat{V}(U) \\
                & \geq \sum_{U\in \bigcup\limits_{T\in \cB}\cB_T} \hat{V}(U) \\
                & \geq \hat{V}^{\DBNR}_{\Rfam_{\hat{V}}}(S),
            \end{align*}
        \end{itemize}
        where the last inequality comes from the fact that 
        $S \subseteq \bigcup\limits_{T\in \cB} T \subseteq \bigcup\limits_{U\in \bigcup\limits_{T\in \cB}\cB_T} U$.
        Hence, $\hat{V}^{\DBNR}_{\Rfam_{\mathrm{op}_{\DBNR}\paren{\hat{V}}}}(S) = \hat{V}^{\DBNR}_{\Rfam_{\hat{V}}}(S)$ and we have
        $\mathrm{op}_{\DBNR}\paren[1]{\hat{V}}=  \mathrm{op}_{\DBNR}\paren[2]{\mathrm{op}_{\DBNR}\paren[1]{\hat{V}}}$.
        
    \subsection{Proof of \autoref{coherent_envelopes}} \label{proof_coherent_envelopes}
        Let $\Rfam = \paren{R_k,\zeta_k}_{k\in \cK}$ be a reference family. 
        We show that $\hat{V}^{\DBNR}_\Rfam$ and $\hat{V}^{\JER}_\Rfam$ are coherent.\\
        This proof relies on Lemma~3 
        from \cite{goeman2021only}: a confidence envelope $\hat V$ is coherent if and 
        only if for all disjoint $S_1, S_2 \subseteq \Nm$, 
        \begin{equation}\label{eq_coherent_lemma}
            \hat V \paren{S_1}\leq \hat V \paren{S_1\cup S_2}\leq \hat V \paren{S_1}+\hat V \paren{S_2}.
        \end{equation}
        Then, with $S_1, S_2 \subseteq \Nm$ disjoint, we have: 
        \begin{itemize}
            \item $\hat{V}^{\DBNR}_\Rfam(S_1) \leq \hat{V}^{\DBNR}_\Rfam(S_1\cup S_2)$ since 
            for all $Q\subseteq \cK$, \[\abs[3]{S_1\cup S_2\setminus\bigcup_{k\in Q}R_k}\geq 
            \abs[3]{S_1\setminus\bigcup_{k\in Q}R_k}.
            \]
            \item $\hat{V}^{\DBNR}_\Rfam(S_1\cup S_2) \leq \hat{V}^{\DBNR}_\Rfam(S_1) + 
            \hat{V}^{\DBNR}_\Rfam(S_2)$ since for all 
            $Q_1, Q_2 \subseteq \cK$, 
            \[\abs[3]{S_1\cup S_2\setminus\bigcup_{k\in Q_1\cup Q_2}R_k}\leq 
            \abs[3]{S_1\setminus\bigcup_{k\in Q_1}R_k} + \abs[3]{S_2\setminus\bigcup_{k\in Q_2}R_k},\]
            and $\sum_{k\in Q_1\cup Q_2}\zeta_k \leq \sum_{k\in Q_1}\zeta_k +\sum_{k\in Q_2}\zeta_k$, 
            in particular with $Q_1$ and $Q_2$ respectively achieving $\hat{V}^{\DBNR}_\Rfam $ 
            for $S_1$ and $S_2$.
        \end{itemize}
        Hence, $\hat{V}^{\DBNR}_\Rfam$ is coherent. Remark that this proof does not use the fact 
        that $S_1$ and $S_2$ are disjoint; the disjoint property of the sets is not required for Lemma~3 from \cite{goeman2021only},
        since if $\hat{V}$ verifies \eqref{eq_coherent_lemma} for all disjoint $S_1,S_2 \subseteq \Nm$, 
        then it also verifies it for all $S_1,S_2 \subseteq \Nm$.\\
        To show that $\hat{V}^{\JER}_\Rfam$ is coherent, we use the same lemma as above with the following inequalities:
        \begin{itemize}
            \item $\hat{V}^{\JER}_\Rfam(S_1) \leq \hat{V}^{\JER}_\Rfam(S_1\cup S_2)$ since 
            for all $A \in \cA(\Rfam)$, $\abs{A\cap S_1}\leq \abs{A\cap (S_1\cup S_2)}$. 
            \item $\hat{V}^{\JER}_\Rfam(S_1\cup S_2) \leq \hat{V}^{\JER}_\Rfam(S_1) + 
            \hat{V}^{\JER}_\Rfam(S_2)$ since for all 
            $A \in \cA(\Rfam)$ such that $\hat{V}^{\JER}_\Rfam(S_1\cup S_2) = \abs{A\cap (S_1\cup S_2)}$, 
            \[\abs{A\cap (S_1\cup S_2)} \leq \abs{A\cap S_1} + \abs{A\cap S_2} \leq \hat{V}^{\JER}_\Rfam(S_1) + \hat{V}^{\JER}_\Rfam(S_2),\]
            where the last inequality comes from the definition of $\hat{V}^{\JER}_\Rfam$.
        \end{itemize}
        Hence, $\hat{V}^{\JER}_\Rfam$ is coherent.

    \subsection{Proof of \autoref{prop_jer_dbnr_extension}} \label{proof_jer_dbnr_extension}
        Let $S\subseteq \Nm$ and $\Rfam = \paren{R_k,\zeta_k}_{k\in \cK}$ a reference family.
        We show that $\tilde{\Rfam} = \paren[1]{\Rfam , \paren{S,|S|}} $ 
        is such that $\hat{V}^{\DBNR}_{\tilde{\Rfam}} = \hat{V}^{\DBNR}_{\Rfam}$ and
        $\hat{V}^{\JER}_{\tilde{\Rfam}}= \hat{V}^{\JER}_{\Rfam}$.\\
        We first show that $\hat{V}^{\DBNR}_{\tilde{\Rfam}} = \hat{V}^{\DBNR}_{\Rfam}$. We have
        $\hat{V}^{\DBNR}_{\tilde{\Rfam}}\leq \hat{V}^{\DBNR}_{\Rfam}$ trivially 
        since $\hat{V}^{\DBNR}_{\tilde{\Rfam}}$ is the minimum of a larger set 
        than $\hat{V}^{\DBNR}_{\Rfam}$.
        Let $Q\subseteq \cK$. Then, we have, for all $T\subseteq \Nm$, 
        \begin{align*}
            \abs[3]{T\setminus\paren[2]{\bigcup_{k\in Q} R_k\cup S}} + \sum_{k\in Q}\zeta_k + |S| 
            &= \abs[3]{T\setminus \bigcup_{k\in Q} R_k} + \sum_{k\in Q}\zeta_k + |S| - \abs[3]{T\setminus\paren[2]{\bigcup_{k\in \tilde{Q}} R_k\cup S^c}}\\
            & \geq \abs[3]{T\setminus\bigcup_{k\in Q} R_k} + \sum_{k\in Q}\zeta_k ,
        \end{align*}
        hence $\hat{V}^{\DBNR}_{\tilde{\Rfam}}\geq \hat{V}^{\DBNR}_{\Rfam}$ and we have $\hat{V}^{\DBNR}_{\tilde{\Rfam}} = \hat{V}^{\DBNR}_{\Rfam}$.\\
        The equality $\hat{V}^{\JER}_{\tilde{\Rfam}}= \hat{V}^{\JER}_{\Rfam}$ is trivial 
        since $\cA(\tilde{\Rfam}) = \cA(\Rfam)$ by definition of $\cA(\cdot)$.
        Hence, reasoning by induction, we have that 
        $\hat{V}^{\DBNR}_{\Rfam_{\hat{V}^{(0)}_\Rfam}} = \hat{V}^{\DBNR}_{\Rfam}$ and
        $\hat{V}^{\JER}_{\Rfam_{\hat{V}^{(0)}_\Rfam}} = \hat{V}^{\JER}_{\Rfam}$ for
        all reference family $\Rfam$.
\section[A coherent confidence envelope that is not a CT procedure]{An example of a coherent confidence envelope that cannot be written as a closed testing procedure} \label{proof_ex_coherent_not_CT}

    Consider \autoref{ex_jer_better}. As shown in \autoref{ex_jer_better}, 
    we have $\hat{V}^{\DBNR}_\Rfam(\Nm) = 2$ and $\hat{V}^{\JER}_\Rfam(\Nm) = 1$.
    We show that the closed testing procedure associated to $\hat{V}^{\DBNR}_\Rfam$ 
    defined in Corollary~1 from \citet{goeman2021only}, that is $\hat{V}^{\CT}_{\mathbf{\varphi}}$ 
    with $\mathbf{\varphi} = \paren{\varphi_S}_{S\subseteq \Nm}$ and $\varphi_S = \ind{\hat{V}^{\DBNR}_\Rfam(S) < |S|}$,
    is different from $\hat{V}^{\DBNR}_\Rfam$. We can easily check that
    for all $S\subseteq \Nm$ such that $|S| \leq 1$, $\varphi_S = 0$ and that for all 
    $S\subseteq \Nm$ such that $|S| \geq 2$, $\varphi_S = 1$. Hence, we have
    \begin{align*}
        \hat{V}^{\CT}_{\mathbf{\varphi}}\paren[1]{\Nm} &= \max_{A\subseteq \Nm}\set{\abs{A}, \inf_{ U\supseteq A} \varphi_U = 0}\\
        &= \max_{A\subseteq \Nm}\set{\abs{A}, |A| \leq 1}\\
        &= 1.
    \end{align*}
    Hence, we have $\hat{V}^{\CT}_{\mathbf{\varphi}}(\Nm) = 1 \neq 2 = \hat{V}^{\DBNR}_\Rfam(\Nm)$, 
    and we have an example of a coherent confidence envelope that cannot be written as a closed testing procedure by the 
    methodology of Corollary~1 from \citet{goeman2021only}, because by \autoref{coherent_envelopes}, for all reference family $\Rfam$,
    $\hat{V}^{\DBNR}_{\Rfam}$ is a coherent confidence envelope. This result obviously shows that Corollary~1 from \citet{goeman2021only} is erroneous.\\
    Furthermore, let us show that $\hat{V}^{\DBNR}_\Rfam$ cannot be written as a closed testing procedure associated to any family of local tests.\\
    Let $\mathbf{\varphi} = \paren{\varphi_S}_{S\subseteq \Nm}$ be a family of local tests and $\hat{V}^{\CT}_\mathbf{\varphi}$ the associated
    closed testing procedure. Assume that $\hat{V}^{\CT}_\mathbf{\varphi} \leq \hat{V}^{\DBNR}_\Rfam$. 
    We show that $\hat{V}^{\CT}_\mathbf{\varphi}(\Nm) < \hat{V}^{\DBNR}_\Rfam(\Nm)$.\\
    We have for all $i\in \set{1,2,3}$, 
    $\hat{V}^{\CT}_\mathbf{\varphi}(R_i) \leq \hat{V}^{\DBNR}_\Rfam(R_i) = \zeta_i = 1$.
    Hence, for all $i\in \set{1,2,3}$, for all $A\subseteq R_i$, $\abs{A} \geq 2$ implies $\inf_{U\supseteq A} \varphi_U = 1$. 
    It implies that for all $A\subseteq \Nm$ such that $\abs{A} \geq 2$, $\inf_{U\supseteq A} \varphi_U = 1$.
    Therefore, we have \[\hat{V}^{\CT}_\mathbf{\varphi}(\Nm) = \max_{A\subseteq \Nm}\set{\abs{A}, \inf_{ U\supseteq A} \varphi_U = 0} \leq 1 < 2 = \hat{V}^{\DBNR}_\Rfam(\Nm).\]
    Hence, $\hat{V}^{\CT}_\mathbf{\varphi} \neq \hat{V}^{\DBNR}_\Rfam$ for all families of local tests $\mathbf{\varphi}$.

\end{document}